\documentstyle{amsart}
\input{bull-art}
\def\descriptionlabel#1{\normalshape\rm\kern-.5em#1}

\newcommand{\mycirc}{\hbox{\raise.25ex\hbox{${\scriptstyle 
\circ}$}}\ }

\renewcommand{\varphi}{\wp}
\newcommand{\R}{\bold{R}} 
\newcommand{\h}{\bold{H}} 
\newcommand{\E}{\bold{E}} 
\newcommand{\p}{{\cal P}}
\newcommand{\m}{{\cal M}}
\newcommand{\BE}{\begin{equation}}
\newcommand{\EE}{\end{equation}}
\newcommand{\BA}{\begin{eqnarray}}
\newcommand{\EA}{\end{eqnarray}}
\newcommand{\BC}[1]{\begin{figure}[ht]\vspace{0.5in}%
\caption #1}
\newcommand{\EC}{\end{figure}}   

\def\0{\kern.5em}		

\newcommand{\leb}[1]{} 

\renewcommand{\Sp}{\bold{S}}

\newtheorem{Theorem}{Theorem}
\newtheorem{Steinitz}{Theorem of Steinitz}

\theoremstyle{remark}
\newtheorem{Notes}{Notes}

\newtheorem{Note}{Note}

\theoremstyle{definition}
\newtheorem{chRast}{Characterization $R^\ast$\defaultfont}

\newtheorem{chR}{Characterization $R$\defaultfont}

\begin{document}
\def\currentvolume{27}
\def\currentissue{2}
\def\currentyear{1992}
\def\currentmonth{October}
\def\copyrightyear{1992}
\def\currentpages{246-251}
\title[Convex hyperbolic and convex polyhedra in the 
sphere]{A 
characterization of convex hyperbolic polyhedra and of
convex polyhedra inscribed in the sphere}
\author[C. D. Hodgson]{Craig~D.~Hodgson}
\address{Mathematics Department, University of
Melbourne, Parkville, Victoria, Australia}
\author{Igor~Rivin}
\address{NEC Research Institute, Princeton, New Jersey 08540
and Mathematics Department, Princeton University,
Princeton, New Jersey 08540}
\author[W. D. Smith]{Warren D. Smith}
\address{NEC Research Institute, Princeton, New Jersey 
08540}
\date{August 30, 1991}
\maketitle

\begin{abstract}
We describe a characterization of convex polyhedra
in $\h^3$ in
terms of their dihedral angles, developed by Rivin. We 
also describe
some geometric and combinatorial consequences of that 
theory. One of
these consequences is a combinatorial characterization of 
convex
polyhedra in $\E^3$ 
all of whose vertices lie on the unit sphere. That 
resolves a
problem posed by Jakob Steiner in 1832.
\end{abstract}

In 1832, Jakob Steiner in his book \cite{st:geom}
asked the following question:

\begin{quote}
In which cases does a convex polyhedron have a 
(combinatorial)
equivalent which is inscribed in, or circumscribed about, 
a sphere?
\end{quote}

This was the 77th of a list of 85 open problems posed by 
Steiner,
of which only numbers 70, 76, and 77 
were still open as of last year.
Apparently Ren\'e Descartes was also interested
in the problem (see \cite{de:solid}). 

Several authors found families of
noninscribable polyhedral types, beginning with Steinitz 
in 1927
(cf. \cite{gru:conp}); all of these families later were 
subsumed by a
theorem of Dillencourt \cite{Dill90b}. In their 1991 book
\cite[problem B18]{CFG}, Croft, Falconer, and Guy
had the following to say:

\begin{quote}
It would of course be nice to characterize the polyhedra of
inscribable type, but as this may be over-optimistic, good 
necessary,
or sufficient, conditions would be of interest.
\end{quote}

Here we announce a full answer to Steiner's question, in 
the sense
that we produce a characterization of inscribable (or
circumscribable) polyhedra that has a number of pleasant 
properties---it
can be checked in polynomial time and it yields a number of
combinatorial corollaries. First we note the following 
well-known
characterization of convex polyhedra proved by Steinitz (cf.
\cite{gru:conp}).

\begin{Steinitz}
A graph is the one-skeleton of a convex
polyhedron in $\E^3$ if and only if it is
a $3$-connected planar graph.
\end{Steinitz}

\begin{Note}
A graph $G$ is {\em $k$-connected} if the complement of any
$k-1$ edges in $G$ is connected.

We will call graphs satisfying the criteria of Steinitz' 
theorem {\em
polyhedral graphs}.

The answer to Steiner's question  stems from the following
characterization of ideal convex polyhedra in hyperbolic 
3-space $\h^3$.
(See \cite{th:gt3m, bear:grps} for the basics of hyperbolic
geometry.) 
\end{Note}
\begin{Theorem}
Let $P$ be a polyhedral graph  with weights $w(e)$ 
assigned to the
edges. Let $P^*$ be the planar dual \RM(or Poincar\'e 
dual\/\RM) 
of $P$, where the edge $e^*$ dual to $e$ is assigned the 
dual weight
$w^*(e^*)= \pi - w(e)$. Then $P$ can be realized as a convex
polyhedron in $\h^3$ with all vertices on the sphere at 
infinity and
with dihedral angle $w(e)$ at every edge $e$ if and only 
if the
following conditions hold\/\rom{:}
\begin{enumerate}
\item{} $0<w^*(e^*)<\pi$ for all edges $e$.
\item{} The sum of dual weights of edges $e_1^*, e_2^*, 
\dots,
e_k^*$ bounding a face in $P^*$ is equal to $2\pi$.
\item{} The sum of dual weights of edges $e_1^*, e_2^*, 
\dots,
e_k^*$ forming a circuit that does not bound a face in 
$P^*$ is strictly
greater than $2\pi$.
\end{enumerate}
\end{Theorem}
\begin{Theorem}
A realization guaranteed by Theorem \RM1 is unique up to 
isometries of
$\h^3$. 
\end{Theorem}

Theorem 1 is proved by Rivin in \cite{ri:ichar}. It
uses the methods of Aleksandrov \cite{Alek50}
and also results and methods developed by Rivin in 
\cite{rivi86,
R1}  and subsequent work. A brief introduction to this 
theory is
given in \S \ref{sec-outline}. A more complete treatment 
is given
in \cite{HR1}. 

\begin{Notes}
Theorems 1 and 5 were recently extended by Rivin to
general hyperbolic polyhedra of finite volume (that is, 
those with
some finite and some ideal vertices).
A characterization of ideal polyhedra with dihedral angles 
not
greater than $\pi/2$ was given by Andreev \cite{Andr70b};
Andreev's result is an easy consequence of Theorem 1.
\end{Notes}

The so-called projective model (or Klein model) of 
hyperbolic 3-space is a
representation of $\h^3$ as the interior of the unit ball 
$B^3$ in the
ordinary Euclidean 3-space $\E^3$. The model has the 
property of being {\it
geodesic}---hyperbolic lines and planes are represented by
Euclidean lines and planes, respectively. Convexity is 
also preserved---a
convex body in $\h^3$ is represented by a convex body in
$B^3$. Thus, hyperbolic convex polyhedra with all vertices 
on the
sphere at infinity correspond precisely to convex 
Euclidean polyhedra
inscribed in the sphere $\Sp^2 = \partial B^3$.
Therefore, a polyhedron
is of inscribable type exactly when it admits an 
edge-weighting that
satisfies the condition in Theorem 1. 

Furthermore (see \cite{gru:conp}), a polyhedron is
inscribable if and only if its planar dual is 
circumscribable, so we
can sum up the characterization as follows. 
\begin{chRast}
A polyhedron $P$ is of circumscribable type if
and only if there exists a weighting $w$ of its edges, 
such that\/\RM:
\begin{enumerate}
\item The weight of any edge satisfies $0<w(e)<1/2.$
\item The total weight of a boundary of a face of $P$ is 
equal to
$1.$
\item The total weight of any circuit not bounding a face is
strictly greater than $1$.
\end{enumerate}
\end{chRast}
\begin{chR}
A polyhedron $P$ is of inscribable type if and only if its 
planar dual
satisfies the conditions \rom{(1)}--\rom{(3)} of 
Characterization $R^*$.
\end{chR}

The following theorem was proved by Smith:
\begin{Theorem}
Given a polyhedral graph $P$, we can decide whether it 
admits a
weighting satisfying Characterization $R^*$ in time 
polynomial in the
number of vertices $N$. More exactly\,\RM:
on an integer Random Access Memory \RM(RAM\/\RM) Machine 
\RM(see \cite{aho74}\RM)
with precision bounded by $O ( \log N )$ bits, the running 
time may be
bounded by $O(N^{5.38})$ operations.
\end{Theorem}

\begin{proof*}{Skeleton of Proof}
Finding the desired weighting 
is a linear program with the number of constraints 
exponential in $N$
and the methods  of \cite{gls:opt} and
\cite{vai:opt} can be used to produce the algorithm of 
Theorem 3. The
algorithm exploits the observation that given a graph with 
{\em
prescribed} weights on the edges, it is possible to 
determine in
polynomal time whether the weights satisfy conditions 
(1)--(3) of
Characterization $R^*$. Given that, a variant of the 
Ellipsoid Method
is seen to yield the desired algorithm. Results of 
\cite{vai:opt}
allow us to improve the asymptotic behavior of the 
algorithm somewhat;
the funny looking exponent $5.38$ stems from the best 
known complexity
result for matrix inversion. 

\begin{Note}[added in proof]
Rivin [20] recently found a much smaller
(linear in $N$) linear program, and hence a simpler 
algorithm.
\end{Note}

Hence, the two realizability questions above may also be 
answered
in polynomial time. For some special classes of graphs,  
it is particularly
easy to decide inscribability. We mention the following 
theorem of
M.~Dillencourt:

\begin{quote}
Any polyhedron whose graph
is 4-connected, is inscribable. 
Also, these graphs are circumscribable.
More graph-theoretic results
can be found in \cite{DS}. 
\end{quote}
\renewcommand{\qed}{}
\end{proof*}

\section{Characterization of hyperbolic polyhedra}
\label{sec-outline}

The work of Aleksandrov \cite{Alek45a, Alek50} gives a
complete characterization of compact convex polyhedra 
in hyperbolic 3-space in terms of the intrinsic hyperbolic 
metric on the
boundary. 
{\em Note\/}:  Aleksandrov's work has now been extended by 
Rivin
\cite{R2} to {\it ideal} convex polyhedra.

Theorem 5 gives an analogous characterization of convex 
hyperbolic  
polyhedra in terms of their dihedral angles. This also 
generalizes the
work of Andreev \cite{Andr70a}. A simple derivation of 
Andreev's
results from Theorem~5 is given by Hodgson in \cite{HR2}.

\subsection{Compact polyhedra}

The material from this section is developed in [18]. See 
[21] for a more detailed
exposition. Consider the {\em Gauss
Map\/} $G$ of a compact convex polyhedron $P$ in Euclidean 
three-dimensional 
space $\E^3$. The map $G$ is a {\sl set-valued function} 
from $P$ to
the unit sphere $\Sp^2$, which assigns to each point $p$ 
the set of
outward unit normals to support planes to $P$ at $p$. 
Thus, the whole
of a face $f$ of $P$ is mapped under $G$ to a single 
point---the
outward unit normal to $f$. An edge $e$ of $P$ is mapped 
to a geodesic
segment $G(e)$ on $\Sp^2$, whose length is easily seen to 
be the exterior
dihedral angle at $e$. A vertex $v$ of $P$ is mapped by 
$G$ to a
spherical polygon $G(v)$, whose sides are the images under 
$G$ of edges
incident to $v$ and whose angles are easily seen to be the 
angles
supplementary to the planar angles of the faces incident 
to $v$; that
is, $G(e_1)$ and $G(e_2)$ meet at angle $\pi-\alpha$ 
whenever $e_1$
and $e_2$ meet at angle $\alpha$. In other words, $G(v)$ 
is exactly the
``spherical polar'' of the link of $v$ in $P$. (The {\em 
link\/} of a vertex is the
intersection of a infinitesimal sphere centered at $v$ 
with $P$, rescaled, so
that the radius is 1.)

Collecting the above observations, it is seen that $G(P)$ is
combinatorially dual to $P$, while metrically it is the
unit sphere $\Sp^2$.

Now apply a similar construction to a convex polyhedron 
$P$ in
$\h^3$. Associate to each vertex $v$ of $P$ a spherical 
polygon $G(v)$
spherically polar to the link of $v$ in $P$. Glue the 
resulting
polygons together into a closed surface, using the rule 
that $G(v_1)$
and $G(v_2)$ are identified isometrically whenever $v_1$ 
and $v_2$ share 
an edge.

The resulting metric space
$G(P)$ is topologically $\Sp^2$ and the complex
is still Poincar\'e dual to $P$. Metrically, however, it 
is no longer
the round sphere. To see this, consider $G(f)$---the 
single common
point of the spherical polygons $G(v_i)$, where $v_i$ is a 
vertex of $f$. 
The angle of $G(v_i)$ incident to $G(f)$ is the exterior 
angle of $f$ at 
$v_i$, and so by the Gauss-Bonnet Theorem, the sum of 
these angles is 
$2\pi+\mathrm{area}(f) \neq 2\pi$. Thus
$G(f)$ is a {\em cone-like singularity\/}, or a {\em cone 
point\/}, with
{\em cone angle\/} greater than $2\pi$. (A cone angle 
equal to $2\pi$
corresponds to a smooth point.)

This analogue of the Gauss map turns out to have rather
remarkable properties. Here is a brief summary:

1.\ The image of a convex Euclidean polyhedron under the 
Gauss
map is always the round sphere $\Sp^2$. In sharp contrast, 
the following
theorem holds.
\begin{Theorem}[Compact Uniqueness]
The metric of
$G(P)$ determines the hyperbolic polyhedron $P$ uniquely 
\RM(up to
congruence\/\RM).
\end{Theorem}

The proof of uniqueness follows the argument used by 
Cauchy in
the proof of his celebrated rigidity theorem for convex 
polyhedra in
$\E^3$ (see [8, 4, or 24]).

2.\ Using the hyperboloid model of hyperbolic 3-space we 
can construct a
model of the map $G$, which is not unlike the well-known 
spherical
polar map. Let $\E^3_1$ denote Minkowski space: $\R^4$ 
equipped with the inner
product of signature $-, +, +, +$. Then $\h^3$ is 
represented by one
sheet of the hyperboloid $\{x \in \E^3_1 \mid \langle x, x 
\rangle = -1\}$, which is the ``sphere of radius 
$\sqrt{-1}\,$'' in $\E^3_1$.
(For a thorough discussion of the hyperboloid model of 
$\h^3$ see
\cite{th:gt3m, bear:grps}.) 

The {\em polar} $P^*$ of a convex polyhedron $P \subset 
\h^3$ consists of
all outward Minkowski unit normals to the support planes 
of $P$.
Each such unit normal vector gives a point in the 
the {\em de Sitter Sphere} $\Sp_1^2= \{x \in \E^3_1 \mid 
\langle x, x
\rangle = 1\}$, which is the ``sphere of radius $1$'' in 
$\E^3_1$.
It turns out that $P^*$ is a convex polyhedron in 
$\Sp_1^2$ and that the
intrinsic metric of $P^*$ is exactly $G(P)$.

\begin{Note}
The de Sitter sphere $\Sp_1^2$ is a semi-Riemannian 
submanifold of
$\E^3_1$ of constant sectional curvature $1$. See 
\cite{onei83} for
further discussion of the geometry of $\E_1^3$ and 
semi-Riemannian
manifolds in general. 
\end{Note}

3.\ We obtain a precise intrinsic characterization of
those surfaces that can arise as $G(P)$ for a compact convex
polyhedron $P$ in $\h^3$. The characterization is quite 
easy to state:

\begin{Theorem}{\defaultfont\bf Characterization Theorem 
for compact
polyhedra.} 
A metric space $(M, g)$ homeomorphic to $\Sp^2$ can arise 
as the Gaussian
image $G(P)$ of a compact convex polyhedron $P$ in $\h^3$
if and only if the following conditions hold\/\RM:

\begin{description}
\item[\rom{(a)}] The metric $g$ has constant curvature $1$ 
away from a
finite collection of cone points $c_i$.
\item[\rom{(b)}] The cone angles at the $c_i$ are greater 
than $2\pi$.
\item[\rom{(c)}] The lengths of closed geodesics of $(M, 
g)$ are all
strictly greater than $2\pi$.
\end{description}
\end{Theorem}

The necessity of (a) and (b) is immediately apparent from 
the
above discussion of $G$. The necessity of (c) is based on
hyperbolic version of Fenchel's theorem (``the total 
geodesic
curvature of a hyperbolic space curve is greater than 
$2\pi$'') and
the ``polarity'' model of the map $G$ sketched in 2. See 
\cite{HR1}
for the details.

The proof of the sufficiency of conditions (a)--(c) is 
based on
Aleksandrov's Invariance of Domain Principle (see 
\cite{Alek39,
Alek50}), which exploits the observation that an open and 
closed
continuous map $f$ from a topological space $A$ into a 
connected
topological space $B$ is necessarily onto. 

Using this idea to prove Theorem 5 requires 
a careful study of the space $\m_n$ of metrics on $\Sp^2$ 
with $n$ cone
points satisfying conditions (a)--(c), of the space $\p_n$ 
of
convex polyhedra in $\h^3$ with $n$ faces, and of the Gauss
map $G : \p_n \to \m_n$.

\subsection{Ideal polyhedra}

The theory of the previous section is extended to noncompact
polyhedra in [22].
Ideal polyhedra can be viewed as ``boundary points'' of 
$\p_n$, and
likewise Theorem 1 can be viewed as a ``limiting case'' of 
Theorem 5.
In particular, a polyhedral graph $P^*$ as in the 
statement of Theorem
1 can be completed to a piecewise-spherical metric on 
$\Sp^2$ by gluing
in a standard round hemi-sphere into each face. It may be
shown that this metric satisfies the conditions (a)--(c) 
of Theorem 5,
{\it except} that it contains closed geodesics of length 
$2\pi$,
corresponding precisely to the equators of the added 
hemi-spheres. 

\begin{Note}
In \cite{R1} the {\em necessity}\ of the conditions of 
Theorem 1 is
established without reference to the characterization of 
compact
polyhedra. 
\end{Note}

\medskip
The techniques used to
prove Theorem 5 are extended to prove Theorem 1 in
\cite{ri:ichar}. The proof involves geometric estimates on 
families of
convex polyhedra in $\h^3$ whose vertices move away to the 
ideal boundary of
$\h^3$ and beyond. The methods actually suffice to produce 
a characterization
of polyhedra of {\it finite volume} in $\h^3$, which 
includes Theorem 1
and Theorem 5 as special cases.
The techniques used to prove
Theorem 4 give only partial uniqueness results for ideal 
polyhedra
(see \cite{R1}). That approach also yields an
algorithm for actually producing an ideal polyhedron in 
$\h^3$ with
prescribed dihedral angles, which runs in time polynomial 
in the number of vertices of
the polyhedron and the number of decimals of accuracy 
required.
In other words this algorithm produces coordinates for a 
convex
inscription of a graph into the unit sphere in $\E^3$. 
This is worthy of note, as the isometric embedding results
of Aleksandrov et al.\ and Theorem 5  
do {\it not} give an effective way to produce a polyhedron
with the desired properties.

\section{Acknowledgments}

The authors would like to thank Brian Bowditch and Mike 
Dillencourt
for helpful discussions. Igor Rivin would like to thank 
Bill Thurston.
He would also like to thank the NEC Research Institute for 
its
hospitality, which made much of this work possible.


\end{document}